# Sophie Germain, mathématicienne extraordinaire: A story stranger than fiction

David Pengelley
Oregon State University
davidp@nmsu.edu

July 22, 2026

**Abstract**

Sophie Germain (1776–1831) was the first woman we know to do significant mathematical research. For two centuries she has been recognized largely for work in elasticity theory. Her research in number theory, on Fermat's Last Theorem, has been known only from a single footnote by Legendre, crediting her with a partial result, today called Sophie Germain's Theorem. However, recent discoveries in her surviving manuscripts and letters have revealed that her accomplishments were much greater, all part of a grand plan she developed and pursued for proving Fermat's Last Theorem in its entirety. We will explain Germain's plan in the context of her life and work.

Imagine a film proposal in which the fictional protagonist, without any apparent formal education or mentoring, and despite strong parental discouragement, teaches herself as a child to university level in mathematics, while a bloody revolution rages for years outside the door. Although our character is prohibited from attending university, as a young woman she engages in gender impersonation to gain access to the world's top researchers, and they provide mentorship encouraging original and substantial research. Along the way she protects the world's top mathematician from an advancing army, fearing he could suffer the same fate as Archimedes, said to have died on the sword of a Roman soldier.

Our film's leading persona becomes the first woman we know of to do important original research in mathematics. She wins a gold medallion prize for her attempted solution to a challenge in elasticity theory set by the academy of sciences. And in number theory she devises a plan for solving another prize problem. This second challenge was an amazing claim about numbers from almost two centuries prior, and it became the most famous unsolved problem in mathematics for another century and a half, until its resolution just a few decades ago. Our protagonist's research carries her plan of solution an impressively long way, as revealed in handwritten manuscripts examined only very recently (just as the imaginary film proposal was being written!). This discovery, although two centuries delayed, considerably raises the stature of the protagonist.

I submit that such a character would be summarily dismissed as utterly implausible by film producers. And yet absolutely everything described above really happened (except not yet the film proposal), and our real life character's name is Sophie Germain[1], whose semiquincentennial we celebrate this year. I will share just a bit of the story [1, 2, 3, 6], and how I fell into recent discoveries through teaching [7, 8, 9, 10].

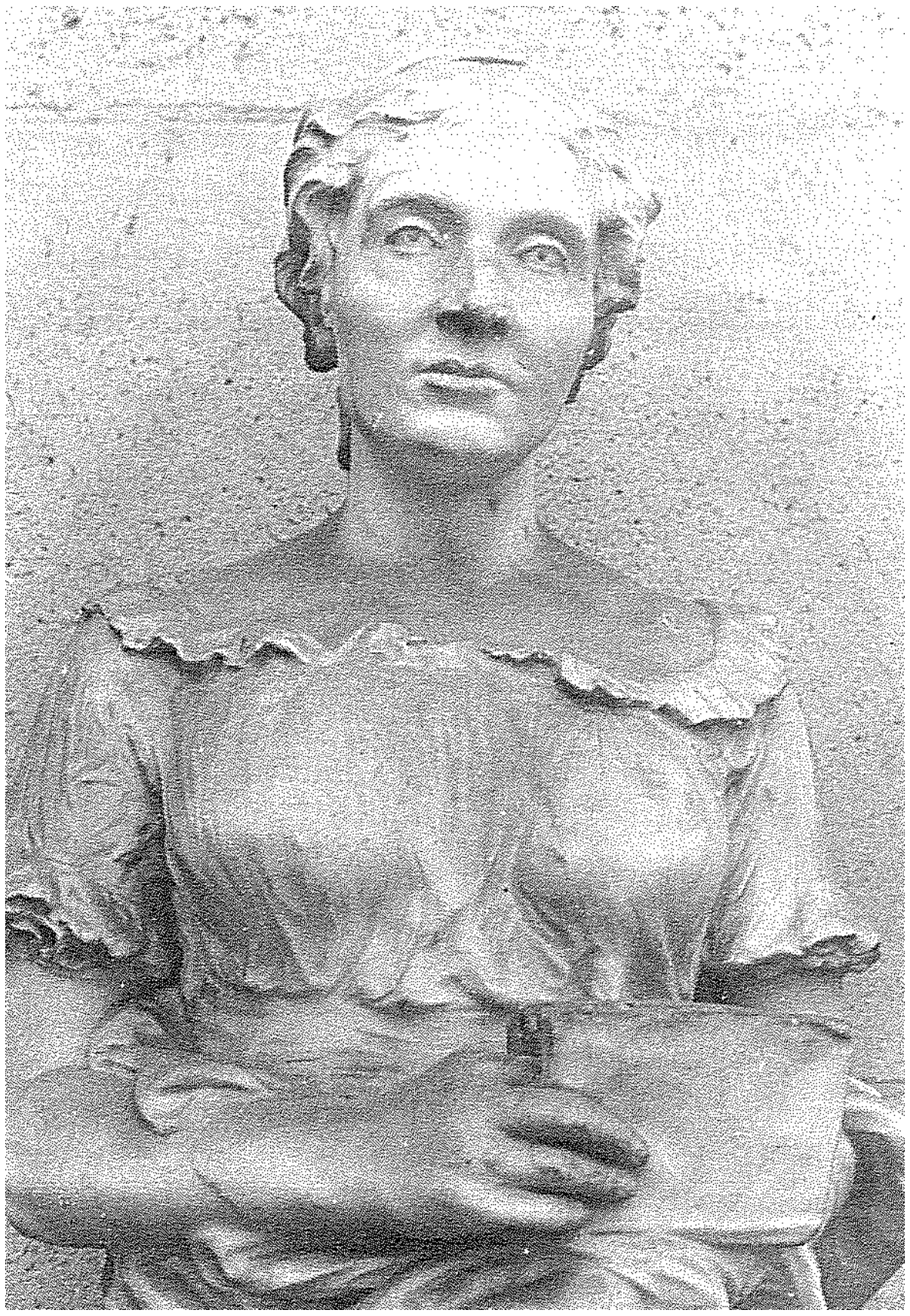

Figure 1: Sophie Germain, a bust by Z. Astruc

Sophie Germain (Figure 1)[2], born in Paris in 1776, was the daughter of a silk merchant. We know only a few remarkable pinpoints about her childhood, mostly from an obituary written by her friend and scientific colleague Guglielmo Libri [1, p. 10ff][3]. At age thirteen she discovered Montucla's *Histoire des mathématiques* in her father's library, as the French Revolution began almost on her doorstep. She continued with a book by Étienne Bézout, and was absorbed in Jacques Cousin's differential calculus during the Reign of Terror in 1793–1794.

Sophie's parents actively prevented her from studying mathematics. According to Libri, she rose to work at night by the light of a lamp, wrapped in blankets, with the ink often frozen in its well, after her parents had removed the fire, her clothes, and candles to force her to stay in bed[3]; but they eventually relented.

[1]Pronounce Ger-main′ as follows (this can be hard for English speakers, but is worth getting right): "G" as in the "t" in "equation", or the "s" in "pleasure" or "vision", or the "g" in "beige"; "er" as in "air"; "main" as in "man", but leave off the "n"!!

[2][1, p. 17]

[3]There is a delightful semi-fictional teenage diary of Sophie Germain [5].

The new École Polytechnique had as instructors some of the world's best researchers. Although as a woman she could not attend, Germain took a most extraordinary step. She obtained their lesson books, focusing particularly on those of Joseph-Louis Lagrange on analysis. Somehow she submitted mathematical comments to Lagrange under the name of Antoine-August LeBlanc, who was an actual student at the École. Lagrange was so impressed that he sought the true author, and went to her to express his astonishment in the most flattering terms [1, p. 11][3].

Sophie Germain became a fervent student of number theory with the publication of Adrien-Marie Legendre's *Essai sur la Théorie des Nombres* in 1798 and of Carl Friedrich Gauss's *Disquisitiones Arithmeticae* in 1801, the two books that founded the modern discipline of number theory. Legendre became a key mentor for Germain.

In 1804 Germain took her next audacious step, writing to Gauss in Germany, again under the pseudonym of the now deceased LeBlanc. She enclosed some of her own research on number theory, particularly on Fermat's Last Theorem (coming below), and Gauss engaged with 'Monsieur LeBlanc' in serious mathematical correspondence. In 1807, Germain intervened with a French general to protect Gauss from Napoleon's army. This action forced her to write to Gauss and admit the true identity of the Monsieur LeBlanc with whom he thought he had been corresponding. Gauss's reply was strikingly different from the negative views of many $19^{\text{th}}$ century scientists and mathematicians on women's abilities.

> 'The taste for the abstract sciences in general and, above all, for the mysteries of numbers, is very rare: this is not surprising, since the charms of this sublime science in all their beauty reveal themselves only to those who have the courage to fathom them. But when a woman, because of her sex, our customs and prejudices, encounters infinitely more obstacles than men, in familiarizing herself with their knotty problems, yet overcomes these fetters and penetrates that which is most hidden, she doubtless has the most noble courage, extraordinary talent, and superior genius.' [1, p. 25]

In 1808, an opportunity presented itself on a track separate from number theory. German physicist E. Chladni created a sensation in Paris demonstrating the intricate nodal patterns arising from vibrations of thin elastic plates, and Napoleon prompted the French Academy of Sciences to set a prize competition for a mathematical explanation. While Germain could attend neither a university nor the Academy, she could enter this competition; submissions were anonymous.

So Germain, despite no formal training, submitted a paper to the competition, claiming to have developed the necessary theory. The Academy twice extended the competition, and in 1816 Germain's third submission received the prize, even though still an incomplete solution. The story is richly told in [1].

The Academy immediately issued a new prize challenge, to solve the open question today known as Fermat's Last Theorem (FLT).

Pierre de Fermat (1601–1665), a century and a half earlier, and Leonhard Euler (1707–1783) in the $18^{\text{th}}$ century, had posed and solved various questions about the natural numbers, setting the stage for Lagrange, Legendre, and Gauss. Fermat made a surprising claim about powers of natural numbers, that no $n^{\text{th}}$ power could be the sum of two other $n^{\text{th}}$ powers, except when $n$ is one or two; formulaically, $z^n = x^n + y^n$ has no solution when $n > 2$. This was particularly astonishing since there are lots of solutions when $n$ is two, such as $5^2 = 3^2 + 4^2$. Fermat's claim came to be known as his 'Last Theorem' because it was the last of his claims to be resolved (in the affirmative), only by applying much sophisticated $20^{\text{th}}$ century mathematics, and developing new techniques of fundamental importance. At the time of the French Academy's new challenge, Fermat's claim had been proven correct for exponents 3 and 4, but was still open for all prime exponents beyond 3. Little did anyone know then that it was one of the hardest questions in mathematics to resolve.

22. Tout se réduit par conséquent à prouver qu'il existe pour chaque valeur de $n$ un nombre premier $\theta$ qui satisfait aux deux conditions mentionnées. Voici un tableau dressé à cet effet pour toutes les valeurs de $n$ moindres que 100.

| $n$ | $\theta$ | $r$ | $n$ | $\theta$ | $r$ |
|---|---|---|---|---|---|
| 3 | 7 = 2. 3+1 | ±1 | 43 | 173 = 4.43+1 | ±(1,80) |
| 5 | 11 = 2. 5+1 | ±1 | 47 | 659 = 14.47+1 | ±(1,12,55,144,249,270,307) |
| 7 | 29 = 4. 7+1 | ±(1,12) | 53 | 107 = 2.53+1 | ±1 |
| 11 | 23 = 2.11+1 | ±1 | 59 | 827 = 14.59+1 | ±(1.20,124,270,337,389,400) |
| 13 | 53 = 4.13+1 | ±(1,23) | 61 | 977 = 16.61+1 | ±(1,52,80,227,252,357,403,439) |
| 17 | 137 = 8.17+1 | ±(1,10,37,41) | 67 | 269 = 4.67+1 | ±(1,82) |
| 19 | 191 = 10.19+1 | ±(1,7,39,49,82) | 71 | 569 = 8.71+1 | ±(1,76,86,277) |
| 23 | 47 = 2.23+1 | ±1 | 73 | 293 = 4.73+1 | ±(1,136) |
| 29 | 59 = 2.29+1 | ±1 | 79 | 317 = 4.79+1 | ±(1,114) |
| 31 | 311 = 10.31+1 | ±(1,6,36,52,95) | 83 | 167 = 2.83+1 | ±1 |
| 37 | 149 = 4.37+1 | ±(1,44) | 89 | 179 = 2.89+1 | ±1 |
| 41 | 83 = 2.41+1 | ±1 | 97 | 389 = 4.97+1 | ±(1,115) |

On voit dans ce tableau que l'équation $r' - r = 1$ n'est satisfaite dans aucun cas, c'est-à-dire qu'il n'y a pas deux restes dont la différence soit égale à l'unité. On voit de même que l'exposant $n$ n'est pas compris parmi les valeurs de $r$. Ainsi la proposition est démontrée, en quelque sorte d'un trait de plume, pour toutes les valeurs de $n$ moindres que 100 (*).

(*) Cette démonstration qu'on trouvera sans doute très ingénieuse, est due à M[lle] *Sophie Germain*, qui cultive avec succès les sciences physiques et mathématiques, comme le prouve le prix qu'elle a remporté à l'Académie sur les vibrations des lames élastiques. On lui doit encore la proposition de l'art. 13 et celle qui concerne la forme particulière des diviseurs premiers de $\alpha$, donnée dans l'art. 11.

Figure 2: Legendre's attribution to Sophie Germain

While Germain never made a submission to the new competition, she worked long and hard at the question. Until recently it has (incorrectly) been assumed that her progress on FLT was exactly as Legendre attributed in a treatise in 1823, the only published record of her work. He first developed a framework that led to a theorem giving a partial resolution of FLT for exponents satisfying certain conditions, and which he explicitly credited to Germain in a footnote (Figure 2)[4]; this was the first general result towards proving Fermat's Last Theorem, and today it is known as Sophie Germain's Theorem (SGT, below). He continued on to some applications of her theorem in certain cases, and he extended some of the results by closer study of key equations, but he did not indicate any role by Germain in the extended results. Recently we have discovered that this is a quite incomplete picture of Germain's accomplishments, which were much greater than what Legendre provides [2, 3].

[4] [3, pp. 645–646]

Germain's approach was built around residues, modular (clock) arithmetic, and congruence. Starting with a fixed natural number $\theta$ called the modulus, any integer $a$ will have a 'residue' (remainder) $r$ upon 'division' by $\theta$, i.e., $a = q\theta + r$ with $q$ an integer and $0 \leq r < \theta$. We call $r$ the residue of $a \bmod \theta$. Assigning to each integer its residue places the integers on a 'clock' with $\theta$ values, and this clock of residues $\bmod \theta$ inherits addition, subtraction, and multiplication. If $a$ and $b$ have the same residue $\bmod \theta$, we write $a \equiv b \pmod{\theta}$ and say that '$a$ is congruent to $b$ modulo $\theta$', which is equivalent to $a - b$ being a multiple of $\theta$. The invention of the notion of congruence is attributed to Gauss, and Sophie Germain was one of the first intensive users of this way of thinking. A key property for Germain's approach is that the clock of residues also admits division if its modulus $\theta$ is prime. For instance, if $\theta = 5$, then the product of the residues 2 and 3 $\bmod \theta$ is 1 (since $2 \cdot 3$ has residue 1, equivalently $2 \cdot 3 \equiv 1 \pmod{5}$). Thus on the $\bmod 5$ clock we can write $1/2 = 3$.

We can now state SGT according to Legendre:

**Sophie Germain's Theorem (SGT).** For an odd prime exponent $p$ in the Fermat equation $z^p = x^p + y^p$, if there exists an auxiliary prime $\theta$ such that there are no two nonzero consecutive $p^{\text{th}}$ power residues modulo $\theta$, nor is $p$ itself a $p^{\text{th}}$ power residue modulo $\theta$, then in any solution to the Fermat equation, one of $x$, $y$, or $z$ must be divisible by $p^2$.

When Legendre proved and credited this theorem to Germain, he also provided a table (Figure 2) that specified a choice of $\theta$ for each $p < 100$ ($p$ is denoted $n$ in his table). Then he listed all the $p^{\text{th}}$ power residues for each $p$ and $\theta$, and verified from these lists that the two hypotheses in SGT are always satisfied. Thus her theorem could successfully be applied for all $p < 100$.

For instance, for $p = 3$, Legendre chooses $\theta = 7$. The nonzero cubic residues $\bmod 7$ are $\{1, 6\}$ (notice that Legendre writes $\pm 1$, which is equivalent), so the two hypotheses are satisfied, and thus according to Germain, in any solution one of $x$, $y$, $z$ must be divisible by $3^2 = 9$.

Naturally, one assumes that Legendre's table reproduces a table of Germain's, i.e., that she applied her theorem as Legendre indicates, by brute force calculation listing $p^{\text{th}}$ powers for chosen $\theta$, to verify the hypotheses for $p < 100$. But it turns out, quite surprisingly, that this is not the case, and the reality is but a small part of a much larger story. I will explain how this discovery arose from teaching.

My students have been learning mathematics directly from studying primary historical sources [7, 8, 9, 10], beginning with an honors course and book cocreated with R. Laubenbacher [4], in which a large theme was FLT, studied via a sequence of primary sources spanning Euclid to Kummer. A highpoint was SGT, the first general result towards proving FLT. But where was the primary source? It seemed that only Legendre's words could speak for Sophie Germain. And yet a close look at Legendre's memoir caused us to wonder whether how he presented SGT was really how Germain would have done. His proof seems to prepare for some strengthened results that he does not attribute to Germain. So what was hers, and what was his? This led us to letters and manuscripts in her own hand that amazingly have survived to the present day, now analyzed in [2, 3].

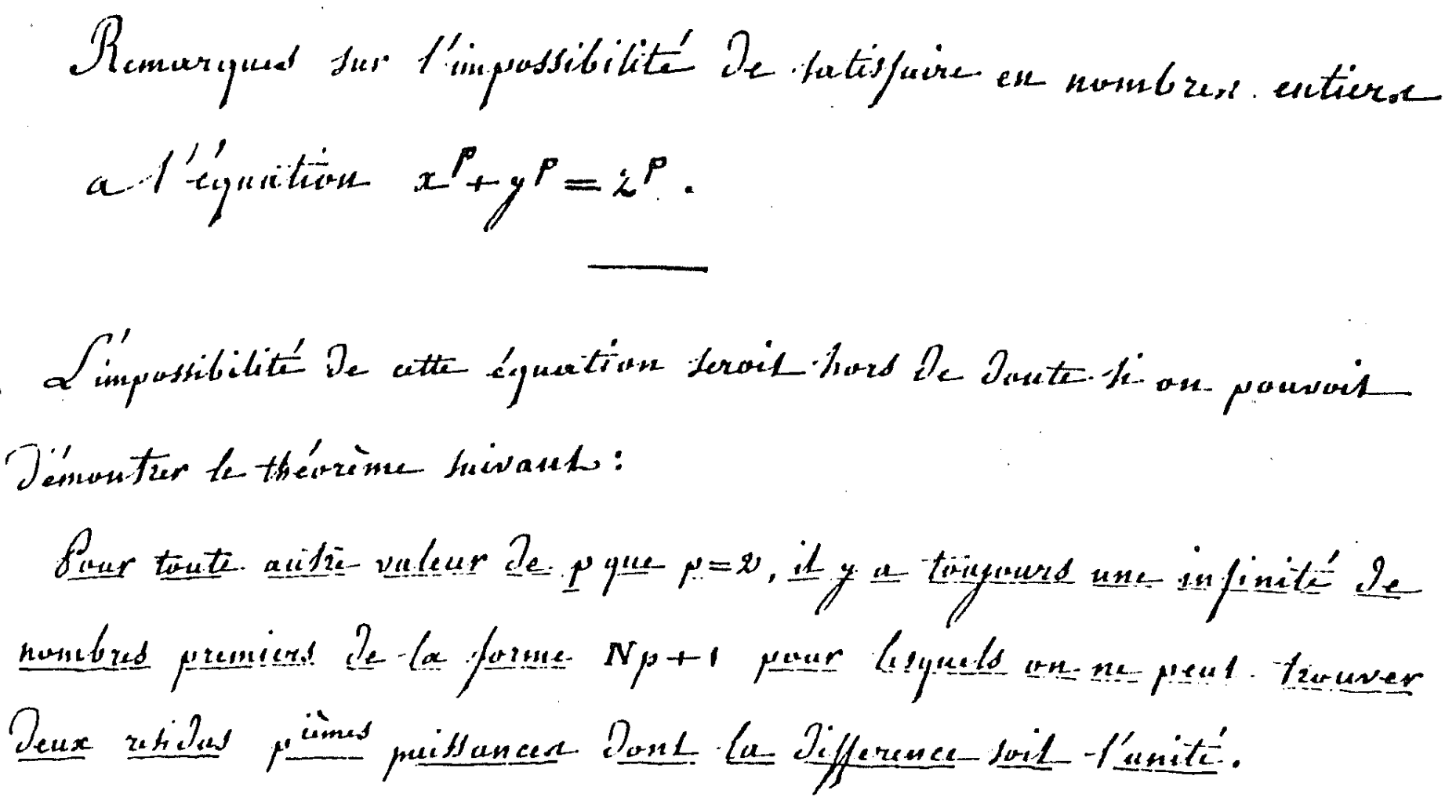

Remarques sur l'impossibilité de satisfaire en nombres entiers a l'équation $x^p + y^p = z^p$.

L'impossibilité de cette équation seroit hors de doute si on pouvoit démontrer le théorème suivant :

Pour toute autre valeur de $p$ que $p = 2$, il y a toujours une infinité de nombres premiers de la forme $Np+1$ pour lesquels on ne peut trouver deux résidus $p^{\text{ièmes}}$ puissances dont la différence soit l'unité.

Figure 3: Sophie Germain's 'Grand Plan' to prove Fermat's Last Theorem

We discovered that Germain's unpublished writings contain not only what Legendre credited to her, but much more. We found that she had what we call her 'Grand Plan' for proving FLT in its entirety, unknown until now, and she carried this plan forward a considerable way. In a nutshell her plan, simple to explain, was as follows (Figure 3)[5].

Fix a prime exponent $p > 2$, and suppose there is a solution $x, y, z$ to the Fermat equation $z^p = x^p + y^p$. Now suppose too that $\theta$ is an auxiliary prime not dividing any of $x, y, z$. Recalling that division is possible mod $\theta$ because $\theta$ is prime, this leads to $\left(\frac{z}{x}\right)^p \equiv 1 + \left(\frac{y}{x}\right)^p \pmod{\theta}$, yielding a nonzero pair of consecutive $p^{\text{th}}$ power residues mod $\theta$. Thus contrapositively, if there are no nonzero consecutive $p^{\text{th}}$ power residues mod $\theta$ (which we will call Condition NC for '**non****c**onsecutivity'; yes, this is the same as one of the conditions in SGT!), then $\theta$ must divide one of $x, y, z$. For instance, for the simple example above where $p = 3$ and $\theta = 7$, Condition NC is true, so one of $x, y, z$ must be divisible by 7.

Sophie Germain's plan was to show that for fixed $p$, there are infinitely many auxiliaries $\theta$ satisfying Condition NC. Then each of these infinitely many $\theta$ must divide one of $x, y, z$, clearly impossible, proving FLT.

Germain's manuscripts show that she developed a web of analyses, conditions, results, methods, and algorithms aiming to prove FLT by establishing the existence of infinitely many auxiliary primes satisfying Condition NC, and that she obtained many partial results. We also found that, unlike Legendre (Figure 2), Germain did not produce a table of brute force calculations to show that SGT can successfully be applied tor all $p < 100$. Rather, for her this follows easily from the greater strength of her theoretical considerations. She applied sophisticated new knowledge and tools, for instance the existence and consequences of primitive roots for any prime modulus, foreshadowed by Euler and first proven by Gauss. And she made extensive use of permutations and their group properties, well before these were formalized.

[5] [3, pp. 649–650, 657–658]

Germain's manuscripts reveal too that she strengthened many of the same results as did Legendre for SGT, but most surprisingly, she and Legendre developed fundamentally different approaches from each other, answering the question we raised from his treatise. Essentially everything he presents beyond SGT, which did not include her grand plan for proving FLT, was also accomplished by her, but independently, without overlap in techniques! This disparate nature of their methods leads us to wonder whether Germain ever received careful reading and feedback from anyone, especially Legendre.

Germain also had a theorem asserting that if the Fermat equation were to have any solutions, they would necessarily be extremely large. For instance, in the case $p = 5$, the theorem would show that one of $x, y, z$ must be at least 39 decimal digits in size! Her proof of this purported theorem contains a flaw, which she later attempted to fix; but this suggests again that perhaps no one else read her work and gave her careful critique. As a woman, Sophie Germain had to work outside the professional academic setting, and we see very concretely here how her work suffered from a lack of collegial interactions comparable to men.

It is perplexing that Legendre does not mention Germain for any of the features beyond SGT that we have uncovered in her work, including her grand plan to prove FLT, her extensions of SGT, and her theoretical techniques and methods. There are several possibilities for explaining this, but the simplest plausible one is that he credited to her the one clearly proven 'theorem' that emerges from all her work; everything else could be considered attempts that didn't culminate in 'results'. Later in the century, other researchers (re)discovered most of what she found beyond SGT, reinventing her methods (not Legendre's), but all while not knowing of her work, since it was never published [3].

In 1819, Sophie Germain wrote a letter to Gauss, after a ten year hiatus in their correspondence (Figure 4)[6]. Her interest in Fermat's Last Theorem, and her dogged tenacity, are illustrated when she writes

> 'Although I have worked for some time on the theory of vibrating surfaces [...], I have never ceased thinking about the theory of numbers. I will give you a sense of my absorption with this area of research by admitting to you that even without any hope of success, I still prefer it to other work which might interest me while I think about it, and which is sure to yield results.' [3, pp. 649–650, 654]

Starting to summarize for Gauss her work towards proving Fermat's Last Theorem she writes 'Voici ce que j'ai trouvé:' ('Here is what I have found:'), and after explaining the essence of her grand plan to him, she addresses the quest for infinitely many auxiliary primes satisfying Condition NC, writing

> 'I have never been able to arrive at the infinity, although I have pushed back the limits quite far by a method of trials too long to describe here.' [3, pp. 649–650, 655]

Finally, she alludes to her 'large size' theorem, in which she at some point discovered a flaw:

[6] [3, pp. 649–650, 654]

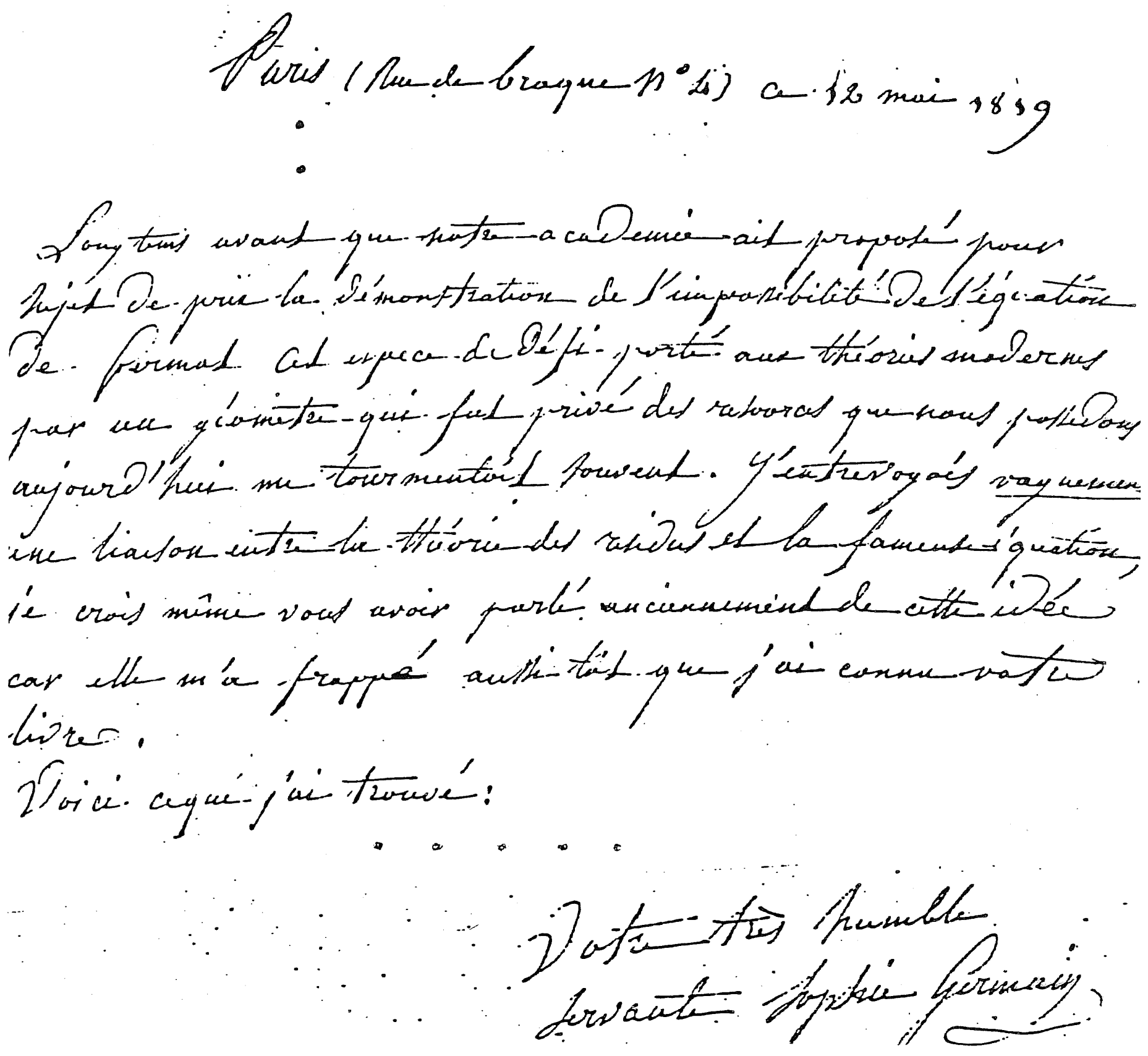
Paris (Rue de braque N° 4) ce 12 mai 1819

Longtems avant que notre academie ait proposé pour sujet de prix la démonstration de l'impossibilité de l'équation de Fermat cet espece de défi porté aux théories modernes par un géomètre qui fut privé des ressources que nous possedons aujourd'hui me tourmentait souvent. J'entrevoyais vaguement une liaison entre la théorie des résidus et la fameuse équation, je crois même vous avoir parlé anciennement de cette idée car elle m'a frappé aussitôt que j'ai connu votre livre.

Voici ce que j'ai trouvé: . . . . .

Votre très humble servante Sophie Germain

Figure 4: Sophie Germain's letter to Gauss, 1819

> 'I have been able to succeed at proving that this equation is not possible except with numbers whose size frightens the imagination.' [3, pp. 649–650, 655]

So could Sophie Germain's grand plan actually work? And did she learn so before her death in 1831? In an undated letter to Legendre she herself actually proves that her approach fails when $p = 3$ (i.e., there are only finitely many auxiliaries satisfying Condition NC) (Figure 5)[7]. And Libri proved this for general $p$ in 1823–1829. Thus while she seemed optimistic about her plan in her letter to Gauss, she likely knew by a few years later that it cannot work. Her experience was thus bittersweet, but this is in the nature of mathematical research, with its thrills and disappointments [3].

To conclude, the incredible and inspiring rediscovered story of Sophie Germain is today linked to the interplay between teaching and history. Teaching with primary sources led us to rich new discoveries about Germain's work, in her own words. This then fed back into teaching, first in a textbook for the honors course that had originally stimulated our search for Germain's words [4, Ch. 4][7], and more recently in an inquiry-based book for a number theory course, built mostly around sources from Sophie Germain herself [8].

A moral and a motto: We can all find gems by seeking and studying primary historical sources. Thank you, Sophie Germain.

[7] [3, pp. 649–650, 665]

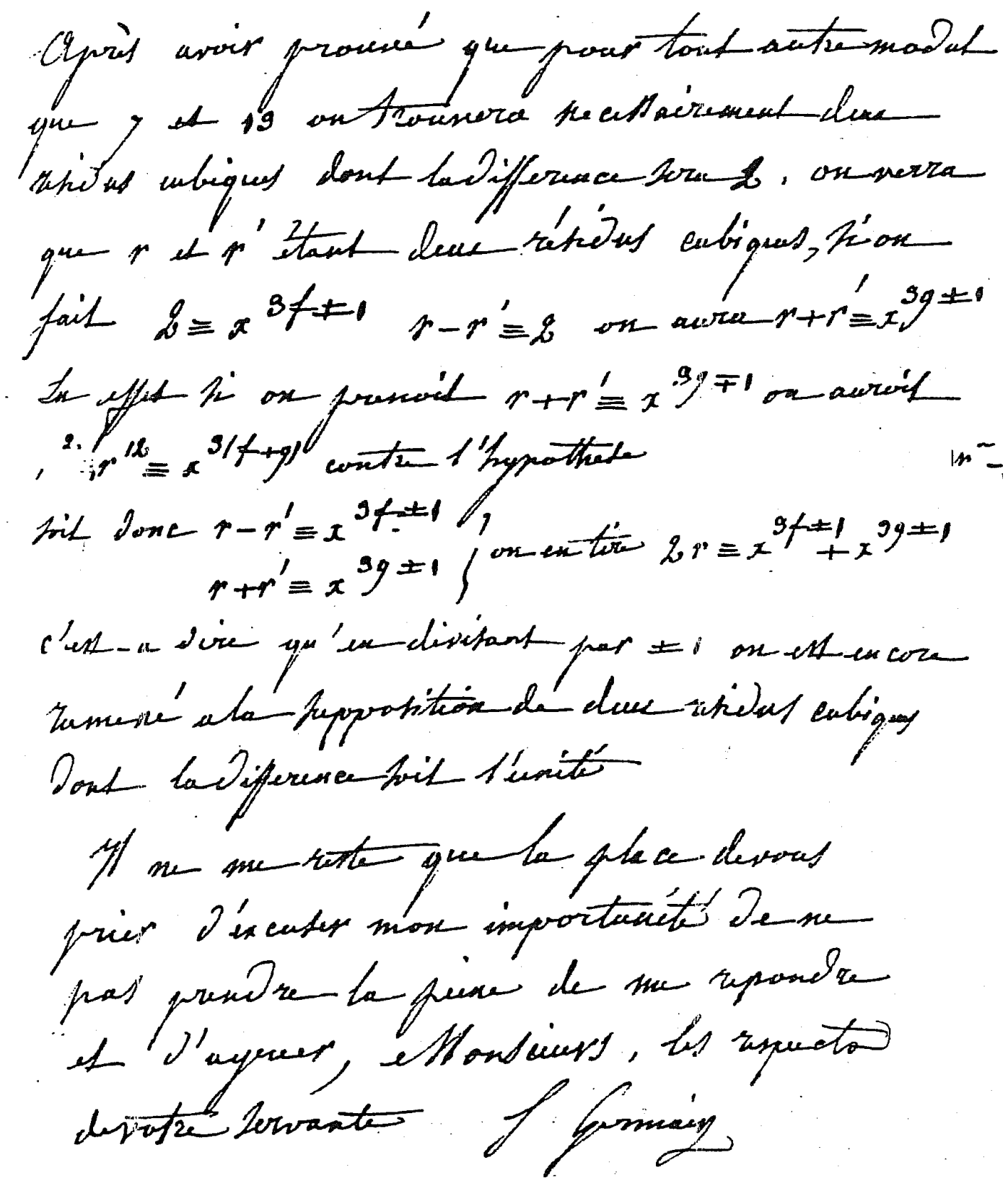

Après avoir prouvé que pour tout autre module que 7 et 13 on trouvera nécessairement deux résidus cubiques dont la différence sera 2, on verra que r et r' étant deux résidus cubiques, si on fait $2 \equiv x^{3f\pm1}$ $r - r' \equiv 2$ on aura $r + r' \equiv x^{3g\pm1}$

En effet si on prenoit $r + r' \equiv x^{3g\mp1}$ on auroit $2 r r' \equiv x^{3(f+g)}$ contre l'hypothèse. Il suit donc $r - r' \equiv x^{3f\pm1}$, $r + r' \equiv x^{3g\pm1}$ on en tire $2r \equiv x^{3f\pm1} + x^{3g\pm1}$ c'est-à-dire qu'en divisant par ±1 on est encore ramené à la supposition de deux résidus cubiques dont la différence soit l'unité.

Il ne me reste que la place de vous prier d'excuser mon importunité de ne pas prendre la peine de me répondre et d'agréer, Monsieurs, les respects de votre servante S Germain

Figure 5: End of Sophie Germain's letter to Legendre

# Sophie Germain, mathématicienne extraordinaire: une histoire plus étrange qu'une fiction

David Pengelley*
Oregon State University
davidp@nmsu.edu

October 27, 2026

**Résumé**

Sophie Germain (1776–1831) a été la première femme que nous connaissions à avoir accompli des recherches mathématiques importantes. Pendant deux siècles, elle a été largement reconnue pour ses travaux sur la théorie de l'élasticité. Ses recherches en théorie des nombres, sur le dernier théorème de Fermat, n'ont été connues que par une seule note de bas de page écrite par Legendre, lui attribuant un résultat partiel, aujourd'hui appelé théorème de Sophie Germain. Cependant, des découvertes récentes dans ses manuscrits et ses lettres conservés ont révélé que ses réalisations étaient bien plus grandes, tout cela faisant partie d'un Grand plan qu'elle a développé et poursuivi pour prouver le dernier théorème de Fermat dans son intégralité. Nous expliquerons le plan de Germain dans le contexte de sa vie et de son travail.

Imaginez une proposition de scénario de film dans laquelle la protagoniste fictive, sans éducation formelle apparente ni mentorat, et malgré une forte dissuasion parentale, s'instruit elle-même en mathématiques depuis l'enfance jusqu'au niveau universitaire, tandis qu'une révolution sanglante fait rage pendant des années devant sa porte. Bien qu'il soit interdit à notre personnage d'aller à l'université, elle s'engage en tant que jeune femme dans une usurpation de genre pour avoir accès aux meilleurs chercheurs du monde, et ils lui offrent un mentorat encourageant des recherches originales et substantielles. En chemin, elle protège le meilleur mathématicien du monde d'une armée qui avance, craignant qu'il ne subisse le même sort qu'Archimède, dont on dit qu'il est mort sous l'épée d'un soldat romain.

Le personnage principal de notre film devient la première femme que nous connaissions à faire d'importantes recherches originales en mathématiques. Elle remporte une médaille d'or pour un prix proposé par l'Académie des sciences, grâce à sa tentative de résolution d'un défi dans la théorie de l'élasticité. Et dans la théorie des nombres, elle conçoit un plan pour résoudre un autre problème faisant l'objet d'un prix académique. Ce deuxième défi était une affirmation étonnante sur les nombres, datant de près de deux siècles auparavant, et il est devenu le problème mathématique non résolu le plus célèbre pendant un autre siècle et demi,

*Je remercie Evelyne Barbin pour une très belle traduction en français.

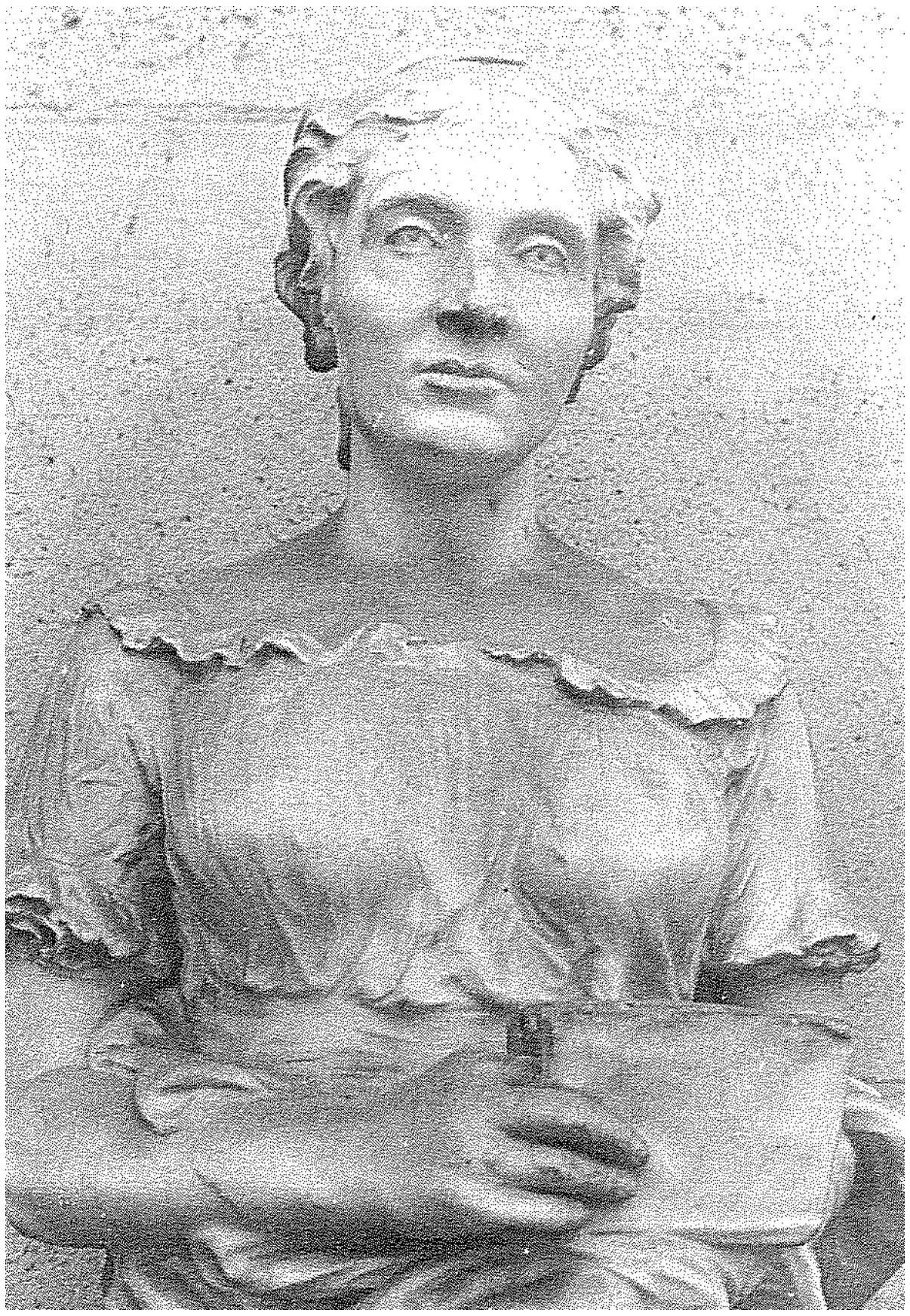

Figure 1: Sophie Germain, le buste par Z. Astruc

jusqu'à sa résolution il y a quelques décennies. La recherche de notre protagoniste expose son plan de résolution d'une manière impressionnante, comme le révèlent des manuscrits de sa main examinés très récemment (juste au moment où la proposition imaginaire du film était en cours d'écriture!). Cette découverte, bien que retardée de deux siècles, accroît considérablement l'importance de la protagoniste.

Je soutiens qu'un tel personnage serait sommairement rejeté comme totalement invraisemblable par les producteurs de films. Et pourtant, absolument tout ce qui a été décrit ci-dessus s'est vraiment produit (sauf, pas encore, la proposition de film), et le nom de notre personnage réel est Sophie Germain, dont nous célébrons le deux cent cinquantième anniversaire cette année. Je vais juste partager ici un peu de cette histoire [1, 2, 3, 6], et dire comment je suis arrivé à des découvertes récentes à travers l'enseignement [7, 8, 9, 10].

Sophie Germain (Figure 1)[1], née à Paris en 1776, était la fille d'un marchand de soie. Nous ne connaissons que quelques points remarquables sur son enfance, principalement d'après une nécrologie écrite par son ami et collègue scientifique Guglielmo Libri [1, p. 10ff][3]. À l'âge de treize ans, elle découvre *l'Histoire des mathématiques* de Montucla dans la bibliothèque de son père, alors que la Révolution française commence, presque sur le pas de sa porte. Elle poursuit avec un livre d'Étienne Bézout, et elle est absorbée par le calcul différentiel de

[1] [1, p. 17]

Jacques Cousin pendant la Terreur, en 1793–1794.

Les parents de Sophie l'empêchaient énergiquement d'étudier les mathématiques. Selon Libri, elle se levait pour travailler la nuit à la lumière d'une lampe, enveloppée dans des couvertures, avec l'encre souvent figée dans son encrier, après que ses parents aient retiré le feu, ses vêtements et ses bougies pour la forcer à rester au lit[2]; mais ils finirent par céder.

La toute nouvelle École polytechnique avait comme professeurs certains des meilleurs chercheurs du monde. Bien qu'en tant que femme, elle ne puisse pas assister aux cours, Germain a entrepris une démarche des plus extraordinaires. Elle a obtenu leurs livres de Leçons, en se concentrant particulièrement sur ceux de Joseph-Louis Lagrange sur l'analyse. D'une manière ou d'une autre, elle a soumis des commentaires mathématiques à Lagrange sous le nom d'Antoine-Auguste Le Blanc, qui était un véritable élève de l'École. Lagrange fut si impressionné qu'il chercha la véritable auteure, et alla la voir pour lui exprimer son étonnement dans les termes les plus flatteurs [1, p. 11][3].

Sophie Germain est devenue une fervente étudiante de la théorie des nombres grâce à la publication de *l'Essai sur la Théorie des Nombres* d'Adrien-Marie Legendre en 1798 et celle des *Disquisitiones Arithmeticae* de Carl Friedrich Gauss en 1801, les deux livres qui ont fondé la discipline moderne de la théorie des nombres. Legendre est devenu un mentor-clé pour Germain.

En 1804, Germain accomplit une démarche audacieuse en écrivant à Gauss en Allemagne, de nouveau sous le pseudonyme de Le Blanc qui est maintenant décédé. Elle a joint une partie de ses propres recherches sur la théorie des nombres, en particulier sur le dernier théorème de Fermat (voir ci-dessous), et Gauss s'est engagé dans une correspondance mathématique sérieuse avec 'Monsieur Le Blanc'. En 1807, Germain intervient auprès d'un Général français pour protéger Gauss de l'armée napoléonienne. Cette action l'oblige à écrire à Gauss et à avouer la véritable identité de Monsieur Le Blanc, avec qui il pensait avoir correspondu. La réponse de Gauss a été remarquablement différente des opinions négatives sur les capacités des femmes de nombreux scientifiques et mathématiciens du XIXe siècle.

> 'Le goût pour les sciences abstraites en général et, surtout, pour les mystères des nombres, est très rare : cela n'a rien de surprenant, puisque les charmes de ces sciences, sublime dans toute leur beauté, ne se révèlent qu'à ceux qui ont le courage de les comprendre. Mais quand une femme, à cause de son sexe, de nos coutumes et préjugés, rencontre infiniment plus d'obstacles que les hommes pour se familiariser avec leurs problèmes épineux, mais surmonte ces entraves et pénètre ce qui est le plus caché, elle a sans aucun doute le courage le plus noble, un talent extraordinaire, et un génie supérieur.' [1, p. 25]

En 1808, une opportunité s'est présentée d'elle-même, sur une voie différente de la théorie des nombres. Le physicien allemand Ernst Chladni a fait sensation à Paris en montrant les lignes nodales complexes résultant de vibrations de fines plaques élastiques, et Napoléon a incité l'Académie des sciences à organiser un prix pour donner une théorie mathématique du phénomène. Bien que Germain ne puisse fréquenter ni une université ni l'Académie, elle pouvait participer à ce Concours; les soumissions étaient anonymes.

[2] Il existe un charmant journal d'adolescente semi-fictif de Sophie Germain [5].

Ainsi, Germain, malgré l'absence de formation formelle, a soumis un texte au Concours, affirmant avoir développé la théorie nécessaire. L'Académie a prolongé le Concours à deux reprises et, en 1816, une troisième soumission de Germain a reçu le prix, même si la solution n'était pas encore complète. L'histoire est richement racontée en [1].

L'Académie a immédiatement lancé un nouveau défi avec un prix pour résoudre la question ouverte, connue aujourd'hui sous le nom du 'dernier théorème de Fermat'.

Pierre de Fermat (1601–1665), un siècle et demi plus tôt, et Leonhard Euler (1707–1783) au XVIIIe siècle, avaient posé et résolu diverses questions sur les nombres naturels, préparant le terrain pour Lagrange, Legendre et Gauss. Fermat a écrit une affirmation surprenante sur les puissances de nombres naturels, selon laquelle aucune puissance $n$ ne peut être la somme de deux autres puissances $n$, sauf lorsque $n$ égale un ou deux; en utilisant une formule, $z^n = x^n + y^n$ n'a pas de solution lorsque $n > 2$. Cela était particulièrement étonnant car il existe de nombreuses solutions lorsque $n$ vaut deux, telle que $5^2 = 3^2 + 4^2$. La demande de Fermat est devenue connue sous le nom de son dernier théorème parce que ceci était la dernière de ses demandes à être résolue (par l'affirmative), seulement en appliquant beaucoup de mathématiques sophistiquées du XXe siècle et en développant de nouvelles techniques d'importance fondamentale. À l'époque du nouveau défi de l'Académie française, l'affirmation de Fermat s'était avérée correcte pour les exposants 3 et 4, mais restait ouverte pour tous les exposants premiers supérieurs à 3. Personne ne savait alors que ceci était l'une des questions les plus difficiles à résoudre en mathématiques.

Bien que Germain n'ait jamais envoyé de soumission pour la nouvelle compétition de l'Académie, elle a travaillé longtemps et dur sur la question. Jusqu'à récemment, il a été supposé (à tort) que ses progrès sur le théorème étaient exactement tels que Legendre les a attribués dans un traité en 1823, ce qui constitue le seul compte rendu publié de son travail. Legendre a d'abord développé un cadre qui a conduit à un théorème donnant une résolution partielle du théorème pour les exposants satisfaisant certaines conditions, et qu'il a explicitement crédité à Germain dans une note en bas de page (Figure 2)[3]; ceci était le premier résultat général pour prouver le dernier théorème de Fermat, et aujourd'hui il est connu sous le nom de théorème de Sophie Germain (ci-dessous). Il a continué avec certaines applications de son théorème dans certains cas, et il a étendu certains des résultats par une étude plus approfondie des équations clés, mais il n'a indiqué aucun rôle de Germain pour les résultats étendus. Récemment, nous avons découvert que ceci est une image assez incomplète des réalisations de Germain, qui étaient beaucoup plus grandes que ce que Legendre fournit [2, 3].

L'approche de Germain était construite à partir des résidus, arithmétique modulaire (horloge), et de la congruence. En commençant avec un nombre naturel fixe $\theta$ appelé module, tout entier $a$ aura un 'résidu' (reste) $r$ lors de la 'division' par $\theta$, c'est-à-dire $a = q\theta + r$ avec $q$ un entier et $0 \leq r < \theta$. Nous appelons $r$ le 'résidu' de $a \bmod \theta$. Attribuer à chaque entier son résidu place les entiers sur une 'horloge' avec $\theta$ valeurs, et cette horloge de résidus $\bmod \theta$ est munie de l'addition, de la soustraction et de la multiplication. Si $a$ et $b$ ont le même résidu $\bmod \theta$, on écrit $a \equiv b \pmod{\theta}$ et on dit que '$a$ est congruent à $b$ $(\bmod \theta)$', ce qui équivaut à ce que $(a - b)$ est un multiple de $\theta$. L'invention de la notion de congruence est attribuée à Gauss, et Sophie Germain fut l'une des premières utilisatrices intensives de cette façon

[3] [3, pp. 645–646]

22. Tout se réduit par conséquent à prouver qu'il existe pour chaque valeur de $n$ un nombre premier $\theta$ qui satisfait aux deux conditions mentionnées. Voici un tableau dressé à cet effet pour toutes les valeurs de $n$ moindres que 100.

| $n$ | $\theta$. | $r$. | $n$. | $\theta$. | $r$. |
|---|---|---|---|---|---|
| 3 | 7 = 2. 3+1 | ±1 | 43 | 173 = 4.43+1 | ±(1,80) |
| 5 | 11 = 2. 5+1 | ±1 | 47 | 659 = 14.47+1 | ±(1,12,55,144,249,270,307) |
| 7 | 29 = 4. 7+1 | ±(1,12) | 53 | 107 = 2.53+1 | ±1 |
| 11 | 23 = 2.11+1 | ±1 | 59 | 827 = 14.59+1 | ±(1.20,124,270,337,389,400) |
| 13 | 53 = 4.13+1 | ±(1,23) | 61 | 977 = 16.61+1 | ±(1,52,80,227,252,357,403,439) |
| 17 | 137 = 8.17+1 | ±(1,10,37,41) | 67 | 269 = 4.67+1 | ±(1,82) |
| 19 | 191 = 10.19+1 | ±(1,7,39,49,82) | 71 | 569 = 8.71+1 | ±(1,76,86,277) |
| 23 | 47 = 2.23+1 | ±1 | 73 | 293 = 4.73+1 | ±(1,136) |
| 29 | 59 = 2.29+1 | ±1 | 79 | 317 = 4.79+1 | ±(1,114) |
| 31 | 311 = 10.31+1 | ±(1,6,36,52,95) | 83 | 167 = 2.83+1 | ±1 |
| 37 | 149 = 4.37+1 | ±(1,44) | 89 | 179 = 2.89+1 | ±1 |
| 41 | 83 = 2.41+1 | ±1 | 97 | 389 = 4.97+1 | ±(1,115) |

On voit dans ce tableau que l'équation $r' - r = 1$ n'est satisfaite dans aucun cas, c'est-à-dire qu'il n'y a pas deux restes dont la différence soit égale à l'unité. On voit de même que l'exposant $n$ n'est pas compris parmi les valeurs de $r$. Ainsi la proposition est démontrée, en quelque sorte d'un trait de plume, pour toutes les valeurs de $n$ moindres que 100 (*).

(*) Cette démonstration qu'on trouvera sans doute très ingénieuse, est due à M[lle] *Sophie Germain*, qui cultive avec succès les sciences physiques et mathématiques, comme le prouve le prix qu'elle a remporté à l'Académie sur les vibrations des lames élastiques. On lui doit encore la proposition de l'art. 13 et celle qui concerne la forme particulière des diviseurs premiers de $\alpha$, donnée dans l'art. 11.

Figure 2: L'attribution de Legendre à Sophie Germain

de penser. Une propriété clé de l'approche de Germain est que l'horloge des résidus admet également la division si son module $\theta$ est premier. Par exemple, si $\theta = 5$, alors le produit des résidus 2 et 3 mod $\theta$ est 1 (puisque $2 \cdot 3$ a pour résidu 1, équivalent à $2 \cdot 3 \equiv 1 \pmod 5$). Ainsi, sur l'horloge mod 5, nous pouvons écrire $1/2 = 3$.

Nous pouvons maintenant énoncer le théorème de Sophie Germain selon Legendre:

**Théorème de Sophie Germain.** Pour un exposant premier impair $p$ dans l'équation de Fermat $z^p = x^p + y^p$, s'il existe un nombre premier auxiliaire $\theta$ tel qu'il n'existe pas deux résidus consécutifs non nuls modulo $\theta$ qui soient puissances $p$-ième, et si $p$ n'est pas lui-même modulo $\theta$ un résidu puissance $p$-ième, alors dans toute solution de l'équation de Fermat, un des $x$, $y$ ou $z$ doit être divisible par $p^2$.

Lorsque Legendre a démontré et crédité ce théorème à Germain, il a également fourni un tableau (Figure 2) qui spécifiait un choix de $\theta$ pour chaque $p < 100$ ($p$ est noté $n$ dans son tableau). Il a ensuite listé tous les résidus de la puissance $p$ pour chaque $p$ et $\theta$; et a vérifié à partir de ces listes que les deux hypothèses du théorème de Sophie Germain sont toujours satisfaites. Ainsi, ce théorème pourrait être appliqué avec succès pour tous les $p < 100$.

Par exemple, pour $p = 3$, Legendre choisit $\theta = 7$. Les résidus cubiques non nuls mod 7 sont $\{1, 6\}$ (notez que Legendre écrit $\pm 1$, ce qui est équivalent), donc les deux hypothèses sont satisfaites, et ainsi, selon Germain, dans toute solution l'une des valeurs $x$, $y$, $z$ doit être divisible par $3^2 = 9$.

Naturellement, on suppose que la table de Legendre reproduit une table de Germain, c'est-à-dire que la mathématicienne a appliqué son théorème comme l'indique Legendre, par un calcul de force brute énumérant les puissances $p$-ième pour le $\theta$ choisi, afin de vérifier les hypothèses pour $p < 100$. Mais il s'avère, de manière assez surprenante, que ce n'est pas le cas, et la réalité n'est qu'une petite partie d'une histoire beaucoup plus grande. Je vais expliquer comment cette découverte est née de l'enseignement.

Mes étudiants ont appris directement les mathématiques en étudiant des sources historiques primaires [7, 8, 9, 10], en commençant par un cours supérieur et un livre conçu avec Reinhard Laubenbacher [4], dans lequel un des grands thèmes était le dernier théorème de Fermat, étudié via une suite de sources primaires allant d'Euclide à Kummer. Un point fort était le théorème de Sophie Germain, qui est le premier résultat général pour prouver le théorème de Fermat. Mais où était la source principale? Il semblait que seules les déclarations de Legendre pouvaient parler pour Sophie Germain. Et pourtant, un examen attentif des mémoires de Legendre nous a conduit à nous demander si la façon dont il a présenté le théorème de Germain était vraiment ce qu'elle avait fait. La preuve de Legendre semble préparer des résultats plus forts, mais qu'il n'attribue pas à Germain. Alors, quelle était celle de Germain, et quelle était la sienne? Cela nous a conduit à des lettres et des manuscrits écrits de sa main qui ont étonnamment survécu jusqu'à nos jours, maintenant analysés dans [2, 3].

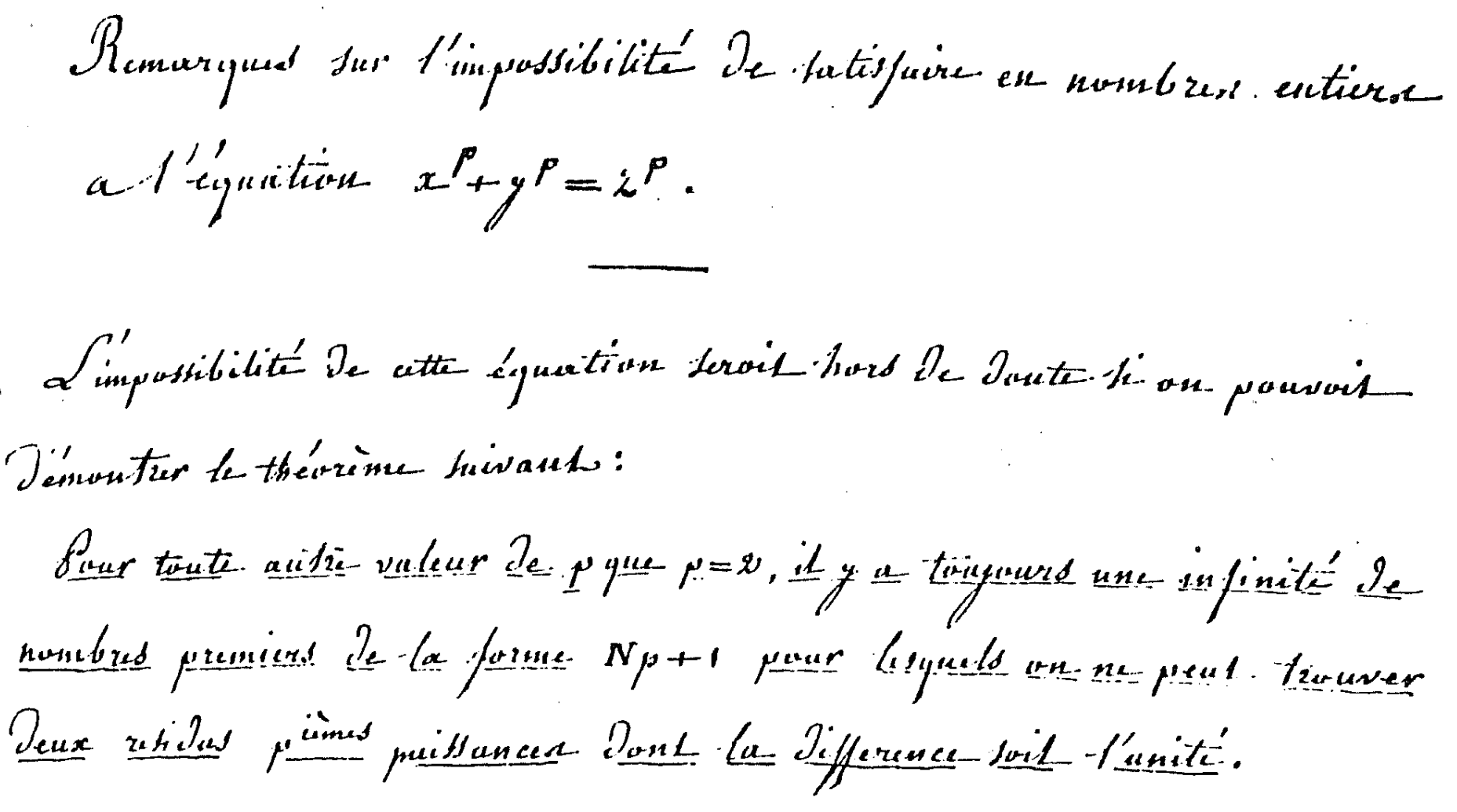
Remarques sur l'impossibilité de satisfaire en nombres entiers à l'équation $x^p + y^p = z^p$.

L'impossibilité de cette équation seroit hors de doute si on pouvoit démontrer le théorème suivant:

Pour toute autre valeur de $p$ que $p = 2$, il y a toujours une infinité de nombres premiers de la forme $Np+1$ pour lesquels on ne peut trouver deux résidus $p^{\text{ièmes}}$ puissances dont la différence soit l'unité.

Figure 3: Le 'Grand plan' de Sophie Germain pour prouver le dernier théorème de Fermat

Nous avons découvert que les écrits inédits de Germain contiennent non seulement ce que Legendre lui créditait, mais bien plus encore. Nous avons constaté qu'elle avait ce que nous appelons son 'Grand plan' pour prouver le dernier théorème de Fermat dans sa totalité, plan inconnu jusqu'à présent, et elle a poursuivi ce plan d'une manière considérable. En résumé, son plan, simple à expliquer, était le suivant (Figure 3)[4].

Posons un exposant premier $p > 2$, et supposons qu'il existe une solution $x, y, z$ à l'équation de Fermat $z^p = x^p + y^p$. Maintenant, supposons aussi que $\theta$ est un nombre premier auxiliaire ne divisant aucun des éléments $x, y, z$. En rappelant que la division est possible $\bmod \theta$ parce que $\theta$ est premier, cela conduit à $\left(\frac{z}{x}\right)^p \equiv 1 + \left(\frac{y}{x}\right)^p \pmod{\theta}$, produisant une paire non nulle de puissances $p$-ième consécutives de résidus $\bmod \theta$. Donc, par contraposition, s'il n'y a pas de puissances $p$-ième consécutives de résidus $\bmod \theta$ non nuls (ce que nous appellerons la condition NC pour '**n**on **c**onsécutivité'; oui, c'est la même que l'une des conditions dans le théorème de Germain!), alors $\theta$ doit diviser l'un des $x, y, z$. Ainsi, pour l'exemple simple ci-dessus où $p = 3$ et $\theta = 7$, la condition NC est vraie, donc l'un des $x, y, z$ doit être divisible par 7.

Le plan de Sophie Germain était de montrer que pour $p$ fixe, il existe une infinité d'auxiliaires $\theta$ satisfaisant la condition NC. Alors chacun de ces infiniment nombreux $\theta$ doit diviser l'un des $x, y, z$, ce qui est clairement impossible, et prouve le dernier théorème de Fermat.

Les manuscrits de Germain montrent qu'elle a développé un réseau d'analyses, de conditions, de résultats, de méthodes et d'algorithmes visant à prouver le dernier théorème de Fermat en établissant l'existence d'une infinité de nombres premiers auxiliaires satisfaisant la condition NC, et qu'elle a obtenu de nombreux résultats partiels. Nous avons également constaté que, contrairement à Legendre (Figure 2), Germain n'a pas produit de tableau de calculs de force brute pour montrer que son théorème peut être appliqué avec succès pour tout $p < 100$. Au contraire, pour elle, cela découle facilement de la plus grande force de ses considérations théoriques. Elle a appliqué des connaissances et des outils sophistiqués, par exemple l'existence de racines primitives pour tout module premier et ses conséquences, préfiguré par Euler et d'abord prouvé par Gauss. Et elle a fait un usage intensif des permutations et de leurs propriétés de groupe, bien avant que celles-ci ne soient formalisées.

Les manuscrits de Germain révèlent aussi qu'elle a renforcé beaucoup des résultats identiques à ceux de Legendre pour son théorème, mais plus étonnamment, elle et Legendre ont développé des approches fondamentalement différentes l'une de l'autre, répondant à la question que nous avions soulevée dans son traité. Essentiellement, tout ce qu'il présente au-delà du théorème de Germain, qui n'incluait pas son Grand plan pour prouver le théorème de Fermat, a également été accompli par elle, mais de manière indépendante, sans chevauchement des techniques! Cette nature disparate de leurs méthodes nous conduit à nous demander si Germain n'a jamais reçu une lecture attentive et des commentaires de quiconque, en particulier de Legendre.

Germain avait également un théorème affirmant que, si l'équation de Fermat devait avoir des solutions, elles seraient nécessairement extrêmement grandes. Par exemple, dans le cas $p = 5$, le théorème montrait que l'un des $x, y, z$ devait avoir une taille d'au moins 39 chiffres! Sa démonstration de ce prétendu théorème contient une faute, qu'elle a ensuite tenté de corriger; mais cela suggère encore que peut-être personne d'autre n'a lu son travail et ne lui a accordé une critique attentive. En tant que femme, Sophie Germain a dû travailler en dehors du cadre académique professionnel, et nous voyons très concrètement ici comment son travail a souffert d'un manque d'interactions collégiales, comparables à celles des hommes.

[4] [3, pp. 649–650, 657–658]

Il est déconcertant que Legendre ne mentionne pas Germain pour aucune des applications qui sont au-delà de son théorème et que nous avons découvertes dans son travail, y compris son Grand plan pour prouver le théorème de Fermat, ses extensions de son théorème et ses techniques et méthodes théoriques. Il y a plusieurs possibilités pour expliquer cela, mais la raison plausible la plus simple est qu'il lui ait crédité le seul 'théorème' clairement prouvé qui émerge de tout son travail; tout le reste pourrait être considéré comme des tentatives qui n'ont pas abouti à des 'résultats'. Plus tard dans le siècle, d'autres chercheurs ont (re)découvert la plupart de ce qu'elle a trouvé au-delà de son théorème, réinventant ses méthodes (pas celles de Legendre), mais tout en ne connaissant pas son travail, puisque celui-ci n'a jamais été publié [3].

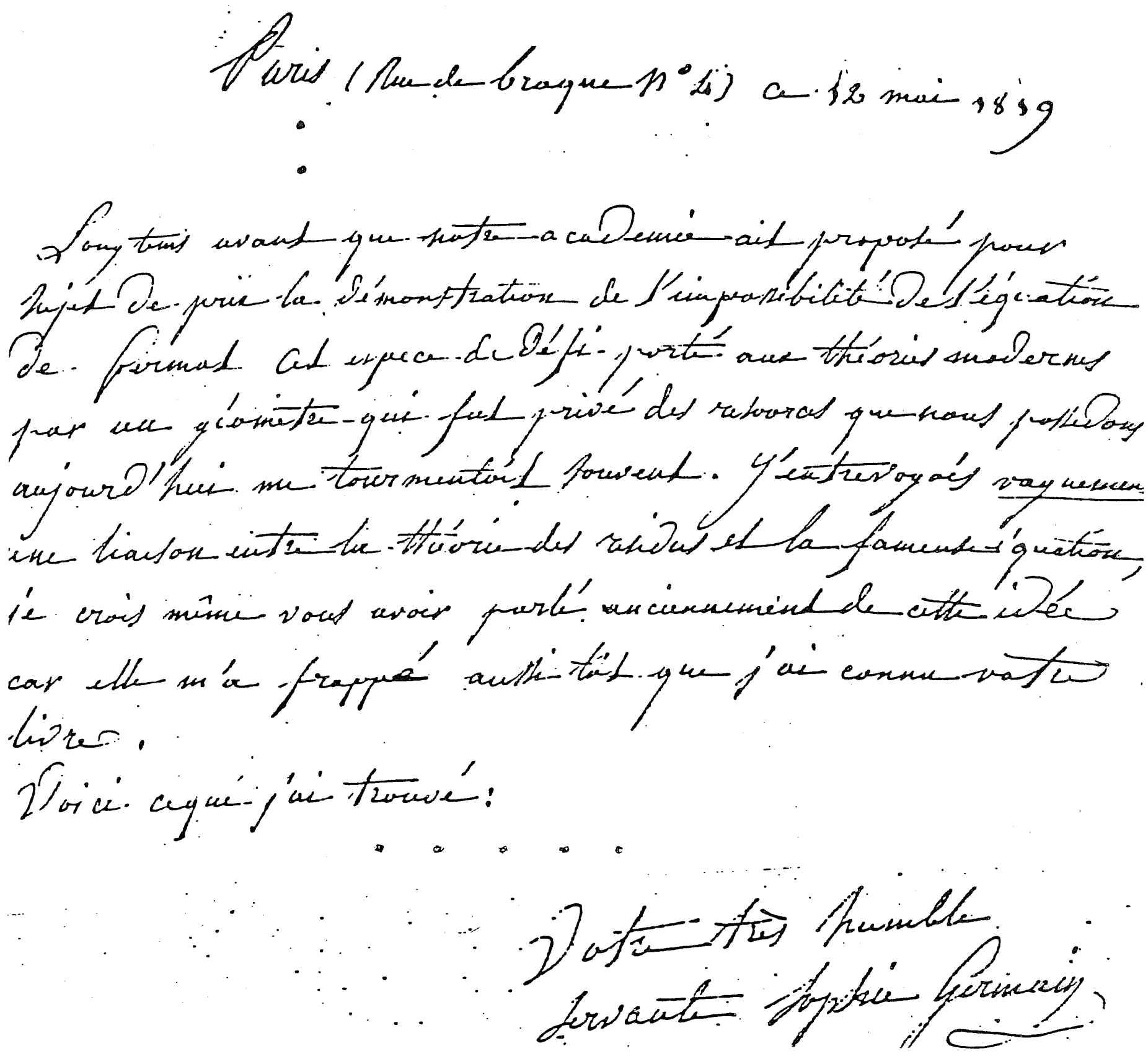
Paris (Rue de braque N° 4) ce 12 mai 1819

Longtems avant que notre academie ait proposé pour sujet de prix la démonstration de l'impossibilité de l'équation de Fermat cet espece de défi porté aux théories modernes par un géomètre qui fut privé des ressources que nous possedons aujourd'hui me tourmentait souvent. J'entrevoyais vaguement une liaison entre la théorie des résidus et la fameuse équation, je crois même vous avoir parlé anciennement de cette idée car elle m'a frappé aussitôt que j'ai connu votre livre.

Voici ce que j'ai trouvé:

. . . . .

Votre très humble servante Sophie Germain

Figure 4: La lettre de Sophie Germain à Gauss (1819)

En 1819, Sophie Germain écrit une lettre à Gauss, après une interruption de dix ans dans leur correspondance (Figure 4)[5]. Son intérêt pour le dernier théorème de Fermat, et sa ténacité obstinée sont illustrés lorsqu'elle écrit:

> 'Quoique j'ai travaillé pendant quelque tems a la théorie des surfaces vibrantes [...], je n'ai jamais cessé de penser a la théorie des nombres. Je vous donnerai une idée de ma préoccupation pour ce genre de recherches en vous avouant que même sans aucune esperance de succès je la prefere a un travail qui me donnerait

[5] [3, pp. 649–650, 654]

necessairement un resultat et qui pourtant m'interresse ... quand j'y pense.' [3, pp. 649–650, 654]

Commençant à résumer pour Gauss son travail en vue de prouver le dernier théorème de Fermat, elle écrit 'Voici ce que j'ai trouvé', et après lui avoir expliqué l'essence de son Grand plan, elle aborde la quête d'une infinité de nombres premiers auxiliaires satisfaisant la condition NC, en écrivant:

> 'Je n'ai jamais pû arriver a l'infini quoique j'ai reculé bien loin les limites par une methode de tatonnement trop longue pour qu'il me soit possible de l'exposer ici.' [3, pp. 649–650, 655]

Enfin, elle fait allusion à son théorème de la 'grande taille', dans lequel elle a découvert une faute à un moment donné:

'... j'ai dû parvenir a prouver que cette équation ne serait possible qu'en nombres dont la grandeur effraye l'imagination.' [3, pp. 649–650, 655]

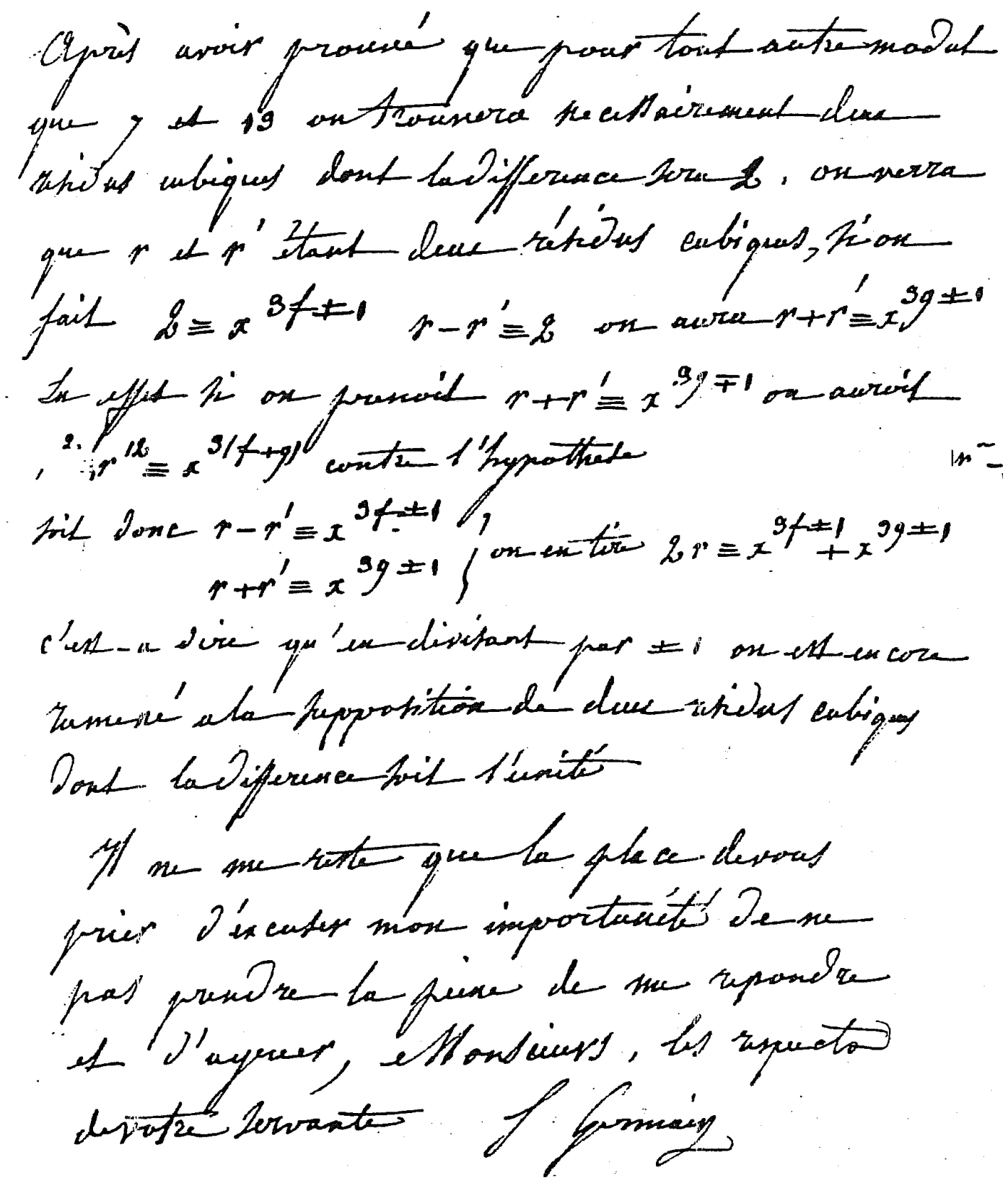
Après avoir prouvé que pour tout autre modul
que 7 et 13 on trouvera necessairement deux
résidus cubiques dont la difference sera 2, on verra
que r et r' étant deux résidus cubiques, si on
fait $2 \equiv x^{3f\pm1}$ $r - r' \equiv 2$ on aura $r + r' \equiv x^{3g\pm1}$
En effet si on prenoit $r + r' \equiv x^{3g\mp1}$ on auroit
$2 r r' \equiv x^{3(f+g)}$ contre l'hypothèse
soit donc $r - r' \equiv x^{3f\pm1}$, $r + r' \equiv x^{3g\pm1}$ on en tire $2r \equiv x^{3f\pm1} + x^{3g\pm1}$
c'est a dire qu'en divisant par ±1 on est encore
ramené a la supposition de deux résidus cubiques
dont la difference soit l'unité
Il ne me reste que la place de vous
prier d'excuser mon importunité de ne
pas prendre la peine de me repondre
et d'agréer, Monsieur, les respects
de votre servante S. Germain

Figure 5: La fin de la lettre de Sophie Germain à Legendre (date non connue)

Ainsi, le Grand plan de Sophie Germain pourrait-il vraiment fonctionner? Et l'a-t-elle appris avant sa mort en 1831? Dans une lettre non datée à Legendre, elle prouve elle-même que son approche échoue lorsque $p = 3$ (c'est-à-dire qu'il n'y a qu'un nombre fini d'auxiliaires satisfaisant la condition NC) (Figure 5)[6]. Et Libri a prouvé ceci pour $p$ quelconque en 1823–1829. Ainsi, bien qu'elle ait semblé optimiste quant à son plan dans sa lettre à Gauss, elle savait probablement quelques années plus tard que cela ne pouvait pas fonctionner. Son expérience était donc douce-amère, mais ceci est dans la nature même de la recherche mathématique, avec ses sensations fortes et ses déceptions [3].

Pour conclure, l'incroyable et inspirante histoire redécouverte de Sophie Germain est aujourd'hui liée à l'interaction entre l'enseignement et l'histoire. Enseigner avec des sources primaires nous a conduit à de nouvelles riches découvertes sur le travail de Germain, selon ses propres mots. Cela a ensuite été réinjecté dans l'enseignement, d'abord dans un manuel pour le cours de spécialisation qui avait initialement stimulé notre recherche des mots de Germain [4, Ch. 4][7], et plus récemment dans un livre d'enquête pour un cours de théorie des nombres, construit principalement à partir de sources de Sophie Germain elle-même [8].

Une morale et une devise: nous pouvons tous trouver des joyaux en recherchant et en étudiant les sources historiques primaires. Merci, Sophie Germain.

# Références

[6] [3, pp. 649–650, 665]